\documentclass[12pt,reqno,a4paper]{amsart}
\usepackage{ae,aecompl}
\usepackage{amssymb,latexsym,pictex,dcpic}
\usepackage[a4paper,top=3cm,right=2.5cm,left=2.5cm,bottom=2.5cm]{geometry}
\usepackage[english]{babel}
\pdfoutput=1

\title{Geometry of the quantum projective plane}
\date{1 December 2009; v2 17 May 2010}

\author[Francesco D'Andrea]{Francesco D'Andrea} 
\address{\flushleft D{\'e}partement de Math{\'e}matique, Universit\'e Catholique de Louvain,
Chemin du Cyclotron 2, B-1348, Louvain-La-Neuve, Belgium}
\email{francesco.dandrea@uclouvain.be}
\author[Giovanni Landi]{Giovanni Landi \\[20pt]}
\address{\flushleft Dipartimento di Matematica e
Informatica, Universit\`{a} di Trieste, Via A. Valerio 12/1, I-34127
Trieste, Italy, and INFN, Sezione di Trieste, Trieste, Italy}
\email{landi@univ.trieste.it}
\thanks{\hspace*{-\parindent}Based on invited talks at the Satellite Conference to the 5th European Congress of Mathematics: ``Noncommutative structures in Mathematics and Physics'', Royal Flemish Academy, Brussels (Belgium), 22-26 July 2008. Published in the Proceedings.}
\newtheorem*{prop}{Proposition}
\newtheorem*{thm}{Theorem}

\newcommand{\dotimes}{\,\dot{\otimes}\,}
\newcommand{\A}{\mathcal{A}}
\newcommand{\E}{\mathcal{E}}
\newcommand{\U}{\mathcal{U}}
\newcommand{\HH}{\mathcal{H}}

\DeclareMathOperator{\SU}{SU}

\newcommand{\N}{\mathbb{N}}
\newcommand{\Z}{\mathbb{Z}}
\newcommand{\R}{\mathbb{R}}
\newcommand{\C}{\mathbb{C}}
\newcommand{\CP}{\mathbb{C}\mathrm{P}}
\newcommand{\Aq}{\mathcal{A}(\mathbb{C}\mathrm{P}^2_q)}
\newcommand{\Oq}{\mathcal{A}(\SU_q(3))}
\newcommand{\Uq}{\mathcal{U}_q(\mathfrak{su}(3))}
\newcommand{\Kq}{\mathcal{U}_q(\mathfrak{u}(2))}

\newcommand{\az}{\triangleright}
\newcommand{\za}{\triangleleft}
\newcommand{\mL}[1]{\mathcal{L}_{#1}}
\newcommand{\wprod}{\wedge_q\mkern-1mu}
\newcommand{\de}{\partial}
\newcommand{\deb}{\bar{\partial}}
\newcommand{\dd}{\mathrm{d}}
\newcommand{\sqbn}[2]{\genfrac{[}{]}{0pt}{}{#1}{#2}}
\newcommand{\inner}[1]{\left<#1\right>}
\newcommand{\nint}{\int\mkern-19mu-\;}
\newcommand{\tr}{\mathrm{Tr}}
\newcommand{\ket}[1]{\left|#1\right>}
\newcommand{\aaz}{\,\textrm{\footnotesize$\blacktriangleright$}\,}

\newcommand{\rPN}{r_N}

\begin{document}

\begin{abstract}
We review some of the geometry of the quantum projective plane 
with emphasis on the construction of a differential calculus and of the Dirac operator (of a spin$^c$-structure). We also report on anti-self-dual connections on line bundles, the spectrum of associated (gauged) Laplacian operators, and on 
classical and quantum characteristic classes.
\end{abstract}

\maketitle

\tableofcontents

\pagebreak

\section*{Introduction}
Among quantum spaces, quantized irreducible flag manifolds occupy a privileged position. On the mathematical side, it is known that the (unique) real covariant differential calculus on such spaces can be realized by commutators with a generalized Dirac operator \cite{Kra04}. 
The simplest of these spaces, the standard Podle\'s sphere, is also a nice toy model for quantum field theory regularization. It turns out that the Dirac operator $D$ has a traceclass resolvent, and as a consequence the basic divergence of the $\phi^4$ theory -- the tadpole diagram --, related to the inverse of the Laplacian, on such a space gives a finite contribution \cite{Oec01}. This is a regularization procedure that does not break the $\SU(2)$ symmetry, but deforms it in a $\SU_q(2)$ symmetry. Aiming at possible applications to physics, we studied a four dimensional example: the quantum projective plane $\CP^2_q$. Further interests for this example come from the fact that it does not admit a spin structure (it is only spin$^c$). We review the construction of the differential calculus on $\CP^2_q$ and of the Dirac operator, the analysis of anti-self-dual connections on line bundles, the spectrum of the associated Laplacians, and the definition of classical and quantum 
characteristic classes.

With $M$ a compact Riemannian spin manifold, the algebra $\A :=C^\infty(M)$ of smooth functions on $M$, the Hilbert space $\HH$ of square integrable spinors with respect to the Riemannian volume form (where functions act by pointwise multiplication), and the Dirac operator $D$ of the Levi-Civita connection, form a triple $(\A ,\HH,D)$ which encodes all the geometrical informations about $M$. It is the prototype of a commutative unital real \emph{spectral triple} \cite{Con95}. As shown in \cite[Theorem 1.2]{Con08} commutative unital real spectral triples are equivalent to oriented compact spin$^c$ manifolds.

In general, a unital spectral triple the datum $(\A ,\HH,D)$ of a (separable) Hilbert space $\HH$, a unital involutive (not necessarily commutative) algebra $\A $ of bounded operators on $\HH$, and a selfadjoint operator $D$ with dense domain in $\HH$. There is a list of properties to be satisfied; in particular one asks that the commutator $[D,f]$ extends to a bounded operator for all $f\in\A $, and that $(D+\mathrm{i})^{-1}$ is a compact operator. The spectral triple is called \emph{even} if $\HH=\HH_+\oplus\HH_-$ is $\Z_2$-graded, the representation of $\A$ is diagonal and the operator $D$ is off-diagonal for this decomposition. The requirement of compact resolvent for the Dirac operator guarantees, for example, that in the even case the twisting of the operator $D^\pm=D|_{\HH_\pm}$ with projections (describing classes in the K-theory of $\A$) are unbounded Fredholm operators: the starting point for the construction of `topological invariants' via index computations \cite{Con94}. Roughly, the bounded commutators condition says that the spectrum of $D$ does not grow too rapidly, while the compact resolvent one says that the specrum of $D$ does not grow too slowly. It is the interplay of the two that imposes stringent restrictions on the geometry and produces spectacular consequences.

Symmetries in noncommutative geometry are encoded in the notion of module algebras (or, dually, comodule algebras), that motivates the definition of a so-called \emph{equivariant} unital spectral triple. Let $\A$ be a complex associative involutive algebra with unity, $(\U,\epsilon,\Delta,S)$ a Hopf $*$-algebra and suppose $\A$ is a left $\U$-module $*$-algebra, which means that there is a left action `$\az$' of $\U$ on $\A$ with properties 
$$
x\az ab=(x_{(1)}\az a)(x_{(2)}\az b)\;,\qquad
x\az 1=\varepsilon(x)1\;,\qquad
x\az a^*=\{S(x)^*\az a\}^*\;,
$$
for all $x\in\U$, $a,b\in\A$. Here we use Sweedler notation for the coproduct, $\Delta(x)=x_{(1)}\otimes x_{(2)}$ with summation understood. The left crossed product $\A\rtimes\U$ is the $*$-algebra generated by $\A$ and $\U$ with crossed commutation relations
$$
xa=(x_{(1)}\az a)x_{(2)}\;,\quad\forall\;x\in\U,\;a\in\A\,.
$$
The data $(\A,\HH,D)$ is called a $\U$-equivariant spectral triple if (i) there is a dense subspace $\mathcal{M}$ of $\HH$ carrying a $*$-representation $\pi$ of $\A\rtimes\U$, (ii) $D$ is a selfadjoint operator with compact resolvent and with domain containing $\mathcal{M}$, (iii) $\pi(a)$ and $[D,\pi(a)]$ extend to bounded operators on $\HH$ for all $a\in\A$, (iv) $[D,\pi(x)]=0$ on $\mathcal{M}$ for any $x\in\U$. In case there is a grading $\gamma$, one further imposes that it commutes with elements of $\U$.

Clearly, any algebra is a module algebra for the trivial action of the bialgebra $\U=\{0\}$. In this case, a $0$-equivariant spectral triple is exactly a usual spectral triple.

In this review, we discuss the geometry of a basic example: the quantum complex projective plane $\CP^2_q$ along the lines of the papers \cite{DDL08b,DL09}. This is defined as a $q$-deformation of the complex projective plane $\CP^2$ seen as the four dimensional real manifold $\SU(3)/\mathrm{U}(2)$.
We shall adopt the notations of \cite{DL09}. Without loss of generality, 
the deformation parameter will be taken to be $0<q<1$.
The symbol $[x]_q$ denotes the $q$-analogue of any $x\in\C$,
$$
[x]_q:=\frac{q^x-q^{-x}}{q-q^{-1}} \;,
$$
and the symbol $[a,b]_q$ denotes the $q$-commutator of two operators $a,b$:
$$
[a,b]_q:=ab-q^{-1}ba \;.
$$
For $n$ a positive integer the $q$-factorial is $[n]_q!:=[n-1]_q![n]_q$, with $[0]_q!:=1$ and  
the $q$-binomial and trinomial are given by
$$
\sqbn{n}{m}_q:=\frac{[n]_q!}{[m]_q![n-m]_q!}\;,
$$
and
$$
[j,k,l]_q!=q^{-(jk+kl+lj)}\frac{[j+k+l]_q!}{[j]_q![k]_q![l]_q!} \;.
$$
Also, the representation symbol $\pi$ will be omitted.

\section{The quantum $\SU(3)$ and $\CP^2$}

Here we present the coordinate algebra of the quantum group $\SU_q(3)$ as an algebra of `true' functions, that is to say as functions on the quantum enveloping algebra $\Uq$: elements in $\A(\SU_q(3))$ are not functions on a group but rather linear maps from $\Uq$ to $\C$. Analogously, we present sections of vector bundles over the quantum projective space $\CP_q^2$ as `true' equivariant maps. This geometric viewpoint allows for a `geometric' description rather that a usual abstract algebraic description. Of course the two approaches are equivalent. 

\subsection{The quantum universal enveloping algebra}
Let $\Uq$ be the compact real form of the Hopf algebra denoted $\breve{U}_q(\mathfrak{sl}(3))$ in Sec.~6.1.2 of~\cite{KS97}.
As a $*$-algebra it is generated by elements $\{K_i,K_i^{-1},E_i,F_i\}_{i=1,2}$,
with $*$-structure $K_i=K_i^*$ and $F_i=E_i^*$, and relations
\begin{gather*}
[K_i,K_j]=0\;,\qquad
[E_i,F_i]=\frac{K_i^2-K_i^{-2}}{q-q^{-1}}\;,  \qquad [E_i,F_j]=0 \quad\mathrm{if}\;i\neq j \;, \\
K_iE_iK_i^{-1}=qE_i\; \qquad K_iE_jK_i^{-1}=q^{-1/2}E_j \quad\mathrm{if}\;i\neq j \;, \\
[E_i,[E_j,E_i]_q]_q=0 \;.
\end{gather*}
It becomes a Hopf $*$-algebra with the following coproduct, counit and antipode:
\begin{gather*}
\Delta(K_i)=K_i\otimes K_i\;,\quad
\Delta(E_i)=E_i\otimes K_i+K_i^{-1}\otimes E_i\;, \quad 
\Delta(F_i)=F_i\otimes K_i+K_i^{-1}\otimes F_i\;, \\
\epsilon(K_i)=1\;,\qquad
\epsilon(E_i)=\epsilon(F_i)=0\;, \\
S(K_i)=K_i^{-1}\;,\qquad
S(E_i)=-qE_i\;,\qquad
S(F_i)=-q^{-1}F_i\;,
\end{gather*}
for $i=1,2$.
For obvious reasons we denote $\U_q(\mathfrak{su}(2))$ the Hopf $*$-subalgebra
of $\Uq$ generated by the elements $\{K_1,K_1^{-1},E_1,F_1\}$, while 
$\Kq$ denotes the Hopf $*$-subalgebra generated by $\U_q(\mathfrak{su}(2))$, $K_1K_2^2$
and $(K_1K_2^2)^{-1}$.

Irreducible representations of $\Uq$ are explicitly described in \cite{DDL08b,DL09}.
Here we simply recall that the representations relevant for our analysis are the
highest weight irreducible representations for which the $K_i$'s are positive operators,
that is to say all irreducible representations that gives representations of $\U(\mathfrak{su}(3))$
when $q\to 1$. They are labelled by a pair $N=(n_1,n_2)\in\N^2$, and basis vectors
are label by a multi-index whose values we indicate with capital letters $I,J$, etc.
In components, $I=(j_1,j_2,m)$ must satisfy suitable constraints, described for example
in \cite{DDL08b}, Eq.~(2.3).

\subsection{Quantized `function algebras'}
The set of linear maps $\Uq\to\C$ is a Hopf $*$-algebra with operations dual to those of $\Uq$.
For $f,g:\Uq\to\C$ we define the product by
$$
(f\cdot g)(x):=f(x_{(1)})g(x_{(2)}) \;,
$$
for all $x\in\Uq$. The unity is the map $1(x):=\epsilon(x)$. Coproduct, counit, antipode and
$*$-involution are given by
\begin{align*}
\Delta(f)(x,y) &:=f(xy) \;,& \hspace{-15mm}
\epsilon(f)&:=f(1) \;,\\
S(f)(x) &:=f(S(x)) \;,& \hspace{-15mm}
f^*(x)&:=\overline{f(S(x)^*)} \;,
\end{align*}
for all $x,y\in\Uq$, and with $\bar{c}$ the complex conjugate of $c\in\C$.

If we focus on maps given by matrix elements of the irreducible representations, we get a Hopf $*$-subalgebra denoted $\Oq$.
A linear basis for $\Oq$ is given by the elements $t^I_J(N)$ defined by
\begin{equation}\label{eq:me}
t^I_J(N)(x):=\rho^N_{I,J}(x) \;,
\end{equation}
where $\rho^N_{I,J}$ is the $(I,J)$-matrix element of the irreducible representation of $\Uq$ with highest weight \mbox{$N=(n_1,n_2)$}.
One can prove that the elements $t^I_J(0,1)$ and $t^I_J(1,0)$ are generators of the algebra (any representation
$\rho^N$ appears as a factor in the 
tensor product $(\rho^{(1,0)})^{n_1}\otimes(\rho^{(0,1)})^{n_2}$ of the
two fundamental representations).

The identification
$$
\begin{pmatrix}
    t^{(0,1,-\frac{1}{2})}_{(0,1,-\frac{1}{2})}(0,1) &
    t^{(0,1,-\frac{1}{2})}_{(0,1,\frac{1}{2})}(0,1) &
    t^{(0,1,-\frac{1}{2})}_{(0,0,0)}(0,1) \\
    t^{(0,1,\frac{1}{2})}_{(0,1,-\frac{1}{2})}(0,1) &
    t^{(0,1,\frac{1}{2})}_{(0,1,\frac{1}{2})}(0,1) &
    t^{(0,1,\frac{1}{2})}_{(0,0,0)}(0,1) \\
    t^{(0,0,0)}_{(0,1,-\frac{1}{2})}(0,1) &
    t^{(0,0,0)}_{(0,1,\frac{1}{2})}(0,1) &
    t^{(0,0,0)}_{(0,0,0)}(0,1)
\end{pmatrix}=
\begin{pmatrix}
    u^1_1 & u^1_2 & u^1_3 \\
    u^2_1 & u^2_2 & u^2_3 \\
    u^3_1 & u^3_2 & u^3_3
\end{pmatrix}
$$
gives an isomorphism between $\Oq$ and the abstract Hopf $*$-algebra generated by elements $u^i_j$
with $i,j=1,2,3$, and defined as follows (cf.~\cite{KS97}, Sec.~9.4). There are relations divided into
commutation relations:
\begin{align*}
u^i_ku^j_k &=qu^j_ku^i_k \;,&
u^k_iu^k_j &=qu^k_ju^k_i \;,&&
\forall\;i<j\;, \\
[u^i_l,u^j_k]&=0 \;,&
[u^i_k,u^j_l]&=(q-q^{-1})u^i_lu^j_k \;,&&
\forall\;i<j,\;k<l\;,
\end{align*}
and the cubic relation
$$
\sum\nolimits_{\pi\in S_3}(-q)^{l(\pi)}u^1_{\pi(1)}u^2_{\pi(2)}u^3_{\pi(3)}=1 \;,
$$
where the sum is over all permutations $\pi$ of the three elements
$\{1,2,3\}$ and $l(\pi)$ is the number of inversions in $\pi$.
The $*$-structure is given by
$$
(u^i_j)^*=(-q)^{j-i}(u^{k_1}_{l_1}u^{k_2}_{l_2}-qu^{k_1}_{l_2}u^{k_2}_{l_1}) \;,
$$
with $\{k_1,k_2\}=\{1,2,3\}\smallsetminus\{i\}$ and
$\{l_1,l_2\}=\{1,2,3\}\smallsetminus\{j\}$, as ordered sets.
Coproduct, counit and antipode are of `matrix' type:
$$
\Delta(u^i_j)=\sum\nolimits_ku^i_k\otimes u^k_j\;,\qquad
\epsilon(u^i_j)=\delta^i_j\;,\qquad
S(u^i_j)=(u^j_i)^*\;,
$$
as expected from \eqref{eq:me}. The elements $t^I_J(N)$ can be written explicitly as
polynomials in the generators $u^i_j$ (cf.~\cite{DL09}, Sec.~2).

The algebra $\Oq$ is a bimodule $*$-algebra for the left and right canonical
actions of $\Uq$. Denoted $\az$ and $\za$ respectively, these actions are the
dual of right (respectively left) multiplication. That is
$$
(x\az f)(y):=f(yx) \qquad\mathrm{and}\qquad (f\za x)(y):=f(xy) \;,
$$
for all $f\in\Oq$ and all $x,y\in\Uq$. Explicitly, on generators:
\begin{align*}
K_i\az u^j_k &=q^{\frac{1}{2}(\delta_{i+1,k}-\delta_{i,k})}u^j_k \;,&
E_i\az u^j_k &=\delta_{i,k} u^j_{i+1}\;, &
F_i\az u^j_k &=\delta_{i+1,k} u^j_i\;, \\
u^j_k\za K_i &=q^{\frac{1}{2}(\delta_{i+1,j}-\delta_{i,j})}u^j_k \;,&
u^j_k\za E_i &=\delta_{i+1,j} u^i_k \;, &
u^j_k\za F_i &=\delta_{i,j} u^{i+1}_k \;.
\end{align*}

The algebra of `functions' on the quantum projective plane $\CP^2_q$ is defined
as the fixed point subalgebra of $\Oq$ for the right action of $\Kq$,
$$
\Aq:=\Oq^{\Kq} \;,
$$
and is a left $\Uq$-module $*$-algebra for the restriction of the left canonical
action. Generators are the elements $p_{ij}:=(u_i^3)^*u_j^3$, and can be arranged
as matrix entries in a projection. For $q=1$, we get a commutative algebra
generated by the matrix entries of a size $3$ and rank $1$ complex projection;
the underlying space is diffeomorphic (as a real manifold) to the projective plane
$\CP^2$ by identifying each line in $\C^3$ with the range of a projection.

\subsection{Hermitian vector bundles}
Roughly speaking, we defined the space $\CP^2_q$ as (the noncommutative analogue of) a quotient
$\SU_q(3)/\mathrm{U}_q(2)$. More generally, one defines vector bundles over $\CP^2_q$
as associated to the principal bundle $\SU_q(3)\to\CP^2_q$ via the representations
of the `structure Hopf algebra' $\Kq$.

It is computationally useful to use the left action $x\mapsto\mL{x}$, of $\Uq$ on $\Oq$,
given by $\mL{x}a:=a\za S^{-1}(x)$; the presence of the antipode yields 
a generalized Leibniz rule:
\begin{equation}\label{glr}
\mL{x}(ab)=(\mL{x_{(2)}}a)(\mL{x_{(1)}}b) \;,
\end{equation}
for $x\in\Uq$ and $a,b\in\Oq$. The pair of commuting actions $\az$
and $\mL{}$ turn $\Aq$ into a left $\Uq\otimes\Uq^{\mathrm{cop}}$-module algebra.
Also, both these left actions are unitary action for the inner product on
$\Oq$ coming from the Haar state $\varphi$:
$$
\varphi\bigl(a^*(x\az b)\bigr)=\varphi\bigl((x^*\az a)^*b\bigr) \qquad\mathrm{and}\qquad
\varphi\bigl(a^*(\mL{x}b)\bigr)=\varphi\bigl((\mL{x^*}a)^*b\bigr) \;,
$$
for all $a,b\in\Oq$ and $h\in\Uq$.

Let $\sigma:\Kq\to\mathrm{End}(\C^n)$ be an $n$-dimensional $*$-representation.
The analogue of (sections of) the vector bundle associated to $\sigma$ is the set
$\E (\sigma)$,
of elements of $\Oq\otimes\C^n$ that are $\Kq$-invariant for the Hopf tensor product
of the actions $\mL{}$ and $\sigma$: 
\begin{multline}\label{eq:Esigma}
\E (\sigma) = \Oq\!\boxtimes_{\sigma}\!\C^n \\ :=
\big\{\varphi \in \Oq\ \otimes \C^n ~\big|~ \big(\mL{h} \otimes \sigma(h)\big)(\varphi) = \epsilon(h) \varphi \,;\;\; \forall\;h\in\Kq\big\} \;.
\end{multline}
As this set is stable under (left and right)
multiplication by an invariant element of $\Oq$, we have that $\E(\sigma)$ is an
$\Aq$-bimodule. For any representation $\sigma$ it can be proved that $\E(\sigma)$
is always projective and finitely generated as one sided (left or right) module~\cite{DL09}.
In addition, since the actions $\mL{}$ and $\az$ commute, it is also a left $\Aq\rtimes\Uq$-module.

Recall that an \emph{Hermitian structure} on a one-sided (say right)
$\A$-module $\E$ is a sesquilinear map $(\,,\,):\E\times\E\to\A$
satisfying $\,(\eta a,\xi b)=a^*(\eta,\xi)b\,$ and $\,(\eta,\eta)\geq 0\,$,
for $\eta,\xi\in\E $ and $a,b\in\A$. 
We also requires the Hermitian structure to be {\it self-dual}, {i.e.} every right $\A$-module homomorphism  $\phi: \E \to \A$ is represented by an element of $\eta \in \E$, by the assignment $\phi(\cdot) = (\eta, \cdot)$, the latter having the correct properties. 

An Hermitian structure exists 
on any finitely generated projective module: if $\E=e\A^k$, with
$e=e^*=e^2$ a size $k$ projection, all Hermitian structure are equivalent
to the one obtained by restricting to $\E$ the standard Hermitian structure
on $\A^k$ given by
\begin{equation}\label{eq:hs}
(\eta,\xi)=\sum\nolimits_{i=1}^k\eta_i^*\xi_i
\end{equation}
for $\eta=(\eta_1,\ldots,\eta_k)$ and $\xi=(\xi_1,\ldots,\xi_k)\in\A^k$.

In particular, for the modules in \eqref{eq:Esigma}, a priori $\eta$ and
$\xi$ have components in $\Oq$, hence $(\eta,\xi)\in\Oq$. Nonetheless,
if $\eta$ and $\xi$ are elements of $\E(\sigma)$, i.e.~they are invariant
under $\mL{}\otimes\sigma$, one checks that $(\eta,\xi)$ is invariant
under the action $\mL{}$,  meaning that $(\eta,\xi)\in\Aq$. Thus, although we
did not give (as yet!) explicitly the isomorphism of $\E(\sigma)$ with a projective
left module, the very same expression \eqref{eq:hs} yields a Hermitian
structure on it. Since $\E(\sigma)$ carries an action of $\Uq$, it is the
$q$-analogue of an equivariant Hermitian vector bundle over $\CP^2$.

A non-degenerate ($\C$-valued) inner product $\inner{\,,\,}$ on $\E(\sigma)$ is obtained
by composing $(\,,\,)$ with the restriction to $\Aq$ of the Haar state:
$$
\inner{\eta,\xi}:=\varphi\big((\eta,\xi)\big) \;,
$$
for $\eta,\xi\in\E $. It is used to define the Hodge $*$-operator, as explained later on.

As $\E (\sigma_1\oplus\sigma_2)\simeq\E(\sigma_1)\oplus\E(\sigma_2)$, it is enough
to focus on irreducible representations of $\Kq$. The highest weight representations
with $K_1$ and $K_2$ positive are classified by a half-integer $\ell$, called the
\emph{spin}, and an integer $N$, called the \emph{charge} (and the dimension over $\C$
is $2\ell+1$). We denote such a representation by the symbol $\sigma_{\ell,N}$ and we
call $\Sigma_{\ell,N}:=\E(\sigma_{\ell,N})$ the associated module.

In particular, we focus on the modules $\Sigma_{0,N}$, with $N\in\Z$, and describe them explicitly
as projective modules. Of course $\Sigma_{0,0}=\Aq$ is the free module with rank $1$.
If $N\neq 0$ we have the following description. Set $z_i:=u_i^3$ (these can be thought of as `coordinates'
on a quantum five-sphere `covering' $\CP^2_q$), and define
\begin{align*}
(\psi_{j,k,l}^N)^*&:=\sqrt{[j,k,l]!}\,z_1^jz_2^kz_3^l\;,
   && \textup{if}\;N>0\; \quad\textup{and with}\quad \; j+k+l=N\,,\\
(\psi_{j,k,l}^N)^*&:=q^{-N+j-l}\sqrt{[j,k,l]!}\,(z_1^jz_2^kz_3^l)^*\;,
   && \textup{if}\;N<0\; \quad\textup{and with}\quad \; j+k+l=-N\,.
\end{align*}
Let $\Psi_N$ be the column vector with components $\psi_{j,k,l}^N$
and $P_N$ the projection given by
\begin{equation}\label{mon-pro}
P_N:=\Psi_N\Psi_N^\dag \;,
\end{equation}
it is of size $\rPN:=\frac{1}{2}(|N|+1)(|N|+2)$. It is shown in~\cite{DL09} that the map
$$
\Sigma_{0,N}\to \Aq^{\rPN}P_{-N}\;,\qquad a\mapsto a\Psi_{-N} \; ,
$$
is an isomorphism of left $\Aq$-modules, while the map
$$
\Sigma_{0,N}\to P_N\Aq^{\rPN}\;,\qquad a\mapsto \Psi_Na \; ,
$$
is an isomorphism of right $\Aq$-modules.

\section{Differential calculus and the Dolbeault-Dirac operator}

Recall that a differential graded algebra (DGA) is the datum $(\Omega^\bullet,\dd)$
of a graded associative algebra $\Omega^\bullet=\bigoplus_{k\geq 0}\Omega^k$ and a
map $\dd:\Omega^\bullet\to\Omega^{\bullet+1}$, which is a cobounday, $\dd^2=0$, and 
a graded derivation:
$$
\dd(\omega_1\omega_2)=(\dd\omega_1)\omega_2+(-1)^{\mathrm{dg}(\omega_1)}\omega_1(\dd\omega_2) \;.
$$
A differential calculus over an associative unital algebra $\A$ is a DGA $(\Omega^\bullet,\dd)$
such that $\Omega^0=\A$ and $\Omega^{k+1}=\mathrm{Span}\{a\dd\omega,\,a\in\A,\,\omega\in\Omega^k\}$
for all $k\geq 0$. For a $^\ast$-calculus (or \emph{real} calculus) there is a graded involution
$^\ast$ on $\Omega^\bullet$ anticommuting with the differential: that is, $^\ast$ is an
involutive antilinear map, $(\omega_1\omega_2)^\ast=(-1)^{\mathrm{dg}(\omega_1)\mathrm{dg}(\omega_2)}
\omega_2^\ast\omega_1^\ast$, and $(\dd\omega)^\ast=-\dd(\omega^\ast)$.

If $\U$ is a Hopf algebra, a real calculus on a $\U$-module $*$-algebra $\A$ is called
$\U$-equivariant if $\Omega^\bullet$ is a left $\U$-module graded $*$-algebra -- with the action respecting the grading --, and $\U$ commuting with the differential.

{}From any spectral triple one can construct a real differential calculus, and the differential
calculus is equivariant if the spectral triple is. In \cite{DDL08b} we proceed in the opposite
way: starting with a (equivariant) differential calculus we construct a Dirac operator and proved 
that the conditions for a (equivariant) unital spectral triple are satisfied.

\subsection{The differential calculus}
Since $\CP^2_q$ is the quantization of a complex manifold, we would like to construct
the analogue of the Dolbeault complex.
Let $\Omega^k:=\bigoplus_{i+j=k}\Omega^{i,j}$ with $\Omega^{i,j}$ the $\ast$-conjugate of $\Omega^{j,i}$,
and assume there is a product such that $\Omega^{\bullet,\bullet}:=\bigoplus_{i,j\geq 0}\Omega^{i,j}$
is a bi-graded $*$-algebra. Given a derivation $\de:\A\to\Omega^{1,0}$, and setting
$\deb a:=-(\de a^\ast)^\ast$, there always exist unique extensions
$\de:\Omega^{\bullet,\bullet}\to\Omega^{\bullet+1,\bullet}$
and $\deb:\Omega^{\bullet,\bullet}\to\Omega^{\bullet,\bullet+1}$ to forms of
arbitrary degree, such that $(\Omega^\bullet,\dd)$ with $\dd:=\de+\deb$, is a real
differential calculus (the relations $\de^2=\de\deb+\deb\de=\deb^2=0$ are equivalent to $\dd^2=0$).
Thus, in order to construct a real differential calculus on $\CP^2_q$, we need (i) two exterior derivations
$\de:\Aq\to\Omega^{0,1}(\CP^2_q)$ and $\deb:\Aq\to\Omega^{0,1}(\CP^2_q)$, (ii) a graded associative $*$-algebra
$\Omega^{\bullet,\bullet}(\CP^2_q)$ such that $\Omega^{0,0}(\CP^2_q)=\Aq$ and such that $\de$ and $\deb$ are conjugated.

In general for $a,b\in\Oq$ and $x\in\Uq$ we have the Leibniz rule in \eqref{glr}.
For $x\in\{E_2,F_2,E_1E_2,F_1F_2\}$, using the $\Kq$-invariance of $\Aq$, from \eqref{glr} we get 
$$
\mL{x}(ab)=(\mL{x}a)b+a(\mL{x}b) \;,\qquad\forall\;a,b\in\Aq\;,
$$
that is, we have four derivations $\Aq\to\Oq$ given by
$\mL{E_2},\mL{F_2},\mL{E_1E_2},\mL{F_1F_2}$; they are exterior derivations since the
image is not in $\Aq$ but rather in $\Oq$. The images of these maps can be identified with
some of the modules previously introduced. {}From now on $\Omega^{0,0}(\CP^2_q):=\Aq$.
Next, we collect the operators in couples. Firstly,
\begin{align*}
\de a  :=
(q^{-\frac{1}{2}}\mL{E_2}a,-q^{\frac{1}{2}}\mL{E_1E_2}a)
\;, \qquad 
\deb a  :=
q^{-\frac{3}{2}}(-q^{-\frac{1}{2}}\mL{F_1F_2}a,q^{\frac{1}{2}}\mL{F_2}a) \;.
\end{align*}
The linear span of the vectors $a\de b$, with $a,b\in\Aq$, is the module $\Sigma_{\frac{1}{2},-\frac{3}{2}}$;
the linear span of the vectors $a\deb b$ is the module $\Sigma_{\frac{1}{2},\frac{3}{2}}$.
The former will be called $\Omega^{1,0}(\CP^2_q)$, the latter $\Omega^{0,1}(\CP^2_q)$.
It is worth mentioning that the $1$-forms $\dd p_{ij}$, with $p_{ij}=u_i^3(u_j^3)^*$
the generators of $\Aq$, are a generating family for $\Omega^1(\CP^2_q)$ as one sided
(left or right) module. For the left module structure the proof is the following.
We have
$$
\de p_{ij}= - q^{-1} (u^3_i)^*\binom{u^2_j}{u^1_j} \;,\qquad
\deb p_{ij}= q^{-1} \binom{-q^{-\frac{1}{2}}(u^1_i)^*}{q^{\frac{1}{2}}(u^2_i)^*}u^3_j \;.
$$
For any $\omega=(v_1,v_2)\oplus (w_1,w_2)\in\Omega^{1,0}(\CP^2_q)\oplus\Omega^{0,1}(\CP^2_q)$ we set
\begin{align*}
a_{ij}(\omega)&:=-q^{1-2j}\bigl\{q^2v_1(u^2_j)^*+ v_2(u^1_j)^*\bigr\}u^3_i \;,\\
b_{ij}(\omega)&:=q^{5-2j}\bigl\{-q^{\frac{1}{2}}w_1(u^3_j)^*u^1_i+q^{-\frac{1}{2}}w_2(u^3_j)^*u^2_i\bigr\} \;.
\end{align*}
These coefficients are right $\Kq$-invariant, i.e.~$a_{ij}(\omega), b_{ij}(\omega)\in\Aq$. Since
$$
\sum\nolimits_jq^{2(a-j)}(u^a_j)^*u^b_j=\sum\nolimits_ju^a_j(u^b_j)^*=\delta_{a,b} \,, \qquad
\sum\nolimits_ju^3_j(u^3_j)^*=\sum\nolimits_jq^{6-2j}(u^3_j)^*u^3_j=1 ,
$$
we get the algebraic identity
$$
\omega= \sum\nolimits_{i,j}\bigl\{a_{ij}(\omega)+b_{ij}(\omega)\bigr\}\dd p_{ij} \;,
$$
which gives the explicit decomposition of any $1$-form $\omega$ in terms of $\dd p_{ij}$.
It goes similarly for the right module structure.

The next point is the definition of the graded $*$-algebra $\Omega^{\bullet,\bullet}(\CP^2_q)$.
Each factor $\Omega^{i,j}(\CP^2_q)$ will be a bimodule associated to a representation of $\Kq$.
We already know which representations correspond to $1$-forms, and for $i+j>0$
the bimodule $\Omega^{i,j}(\CP^2_q)$ must be isomorphic to the tensor product over $\Aq$
of $i$ copies of $\Omega^{1,0}(\CP^2_q)$ and $j$ copies of $\Omega^{0,1}(\CP^2_q)$. Since
$\E(\sigma_1)\otimes_{\Aq}\E(\sigma_2)\subset\E(\sigma_1\otimes\sigma_2)$
for any pair of representations $\sigma_1,\sigma_2$, it is not difficult
to guess the representations corresponding to forms with degree
higher than $1$. They are listed in Table \ref{tab}.

\begin{table}[t]

\begin{center}
\begin{tabular}{ccc}
\begindc{\commdiag}[30]
 \obj(3,5)[A]{$\;\;\;\Omega^{0,0}$}
 \obj(2,4)[B]{$\;\;\;\Omega^{0,1}$}
 \obj(4,4)[C]{$\;\;\;\Omega^{1,0}$}
 \obj(1,3)[D]{$\;\;\;\Omega^{0,2}$}
 \obj(3,3)[E]{$\;\;\;\Omega^{1,1}$}
 \obj(5,3)[F]{$\;\;\;\Omega^{2,0}$}
 \obj(2,2)[G]{$\;\;\;\Omega^{1,2}$}
 \obj(4,2)[H]{$\;\;\;\Omega^{2,1}$}
 \obj(3,1)[I]{$\;\;\;\Omega^{2,2}$}
 \mor{A}{B}{}
 \mor{A}{C}{}
 \mor{B}{D}{}
 \mor{B}{E}{}
 \mor{C}{E}{}
 \mor{C}{F}{}
 \mor{D}{G}{}
 \mor{F}{H}{}
 \mor{E}{G}{}
 \mor{E}{H}{}
 \mor{G}{I}{}
 \mor{H}{I}{}
\enddc
& \raisebox{0.84in}{$=\!\!$} &
\begindc{\commdiag}[30]
 \obj(3,5)[A]{$(0,0)$}
 \obj(2,4)[B]{$(\frac{1}{2},\frac{3}{2})$}
 \obj(4,4)[C]{$(\frac{1}{2},-\frac{3}{2})$}
 \obj(1,3)[D]{$(0,3)$}
 \obj(3,3)[E]{$(1,0)\oplus (0,0)$}
 \obj(5,3)[F]{$(0,-3)$}
 \obj(2,2)[G]{$(\frac{1}{2},\frac{3}{2})$}
 \obj(4,2)[H]{$(\frac{1}{2},-\frac{3}{2})$}
 \obj(3,1)[I]{\smash[b]{$(0,0)$}}
 \mor{A}{B}{}
 \mor{A}{C}{}
 \mor{B}{D}{}
 \mor{B}{E}{}
 \mor{C}{E}{}
 \mor{C}{F}{}
 \mor{D}{G}{}
 \mor{F}{H}{}
 \mor{E}{G}{}
 \mor{E}{H}{}
 \mor{G}{I}{}
 \mor{H}{I}{}
\enddc
\end{tabular}
\end{center}

\caption{In the diamond on the right, in position $(i,j)$ we put
spin and charge $(\ell,N)$ of the representation corresponding to the
bimodule $\Omega^{i,j}$.}\label{tab}
\end{table}
 
As a way of illustration of the general strategy, we give here explicitly the product
of antiholomophic forms (the part that is relevant in the construction
of the Dirac operator). The complete construction can be found in \cite{DL09}.
We have $\Omega^{0,2}(\CP^2_q)=\Sigma_{0,3}$, and an antiholomorphic form will
be a triple $\omega=(a,v,b)$, where $a\in\Omega^{0,0}(\CP^2_q)$ and $b\in\Omega^{0,2}(\CP^2_q)$
are `scalars' and $v=(v_+,v_-)\in\Omega^{0,1}(\CP^2_q)$ is two component vector.
The fact that the product is graded, and that multiplication by scalars
is given by the multiplication in $\Oq$, restrict the possible products
to the following form
$$
\omega\cdot\omega'=(aa',av'+va', ab'+ba'+v\wprod v') \;,
$$
where $\wprod:\Omega^{0,1}(\CP^2_q)\otimes\Omega^{0,1}(\CP^2_q)\to\Omega^{0,2}(\CP^2_q)$ is a linear
map one needs to determine. This can be constructed by using the
Clebsch-Gordan decomposition rules for the product of two spin $\frac{1}{2}$
representations of $\U_q(\mathfrak{su}(2))$. Modulo a normalization,
that we fix according to the notations in \cite{DDL08b}, there is only
one such a product: 
$$
v\wprod v':=\tfrac{2}{[2]_q}(q^{\frac{1}{2}}v_+v_-'-q^{-\frac{1}{2}}v_-v_+') \;.
$$
Notice that equivariance of the differential calculus is straightforward if
we make $\Uq$ to act on forms by lifting diagonally the left action $\az$ on functions.

To get an involution we use the fact that the spin $1/2$ (resp.~spin $1$)
representation of $\U_q(\mathfrak{su}(2))$ is quaternionic (resp.~real).
Rephrased in terms of the representations $\sigma_{\frac{1}{2},N}$ and
$\sigma_{1,N}$ of $\Kq$ we have the following proposition, which takes into
account the fact that real/quaternionic structures change sign to the label $N$.

\begin{prop}\label{lemma:J}
Let $V_{\ell,N}=\C^{2\ell+1}$ be the vector space which underlies the representation $\sigma_{\ell,N}$ of $\Kq$.
An antilinear map $J:V_{\ell,N}\to V_{\ell,-N}$ satisfying $J^2=(-1)^{2\ell}$ and such that 
\begin{equation}\label{eq:J}
J\sigma_{\ell,N}(h)=\sigma_{\ell,-N}(S(h)^*)J, 
\end{equation}
for any $h\in\Kq$, is given, for $\ell=0,\frac{1}{2},1$, by
$$
Ja=a^* \;,\quad
J(v_1,v_2)=(-q^{-\frac{1}{2}}v_2^*,q^{\frac{1}{2}}v_1^*) \;,\quad
J(w_1,w_2,w_3)=(-q^{-1}w_3^*,w_2^*,-qw_1^*) \;,
$$
for any $a\in V_{0,N}$, $v\in V_{\frac{1}{2},N}$ and $w\in V_{1,N}$ respectively.
\end{prop}

\noindent
The operator $J$ is extended to $\omega\in\Omega^{\bullet,\bullet}(\CP^2_q)$ by
composing the omonymous map $J$ on the vectorial part of $\omega$
with the $*$-involution on $\Oq$ (from \eqref{eq:J} it follows that $J$
maps forms to forms).
A graded involution $^\star$ on $\Omega^{\bullet,\bullet}(\CP^2_q)$ is then given by the map
$$
(\omega^\star)_{i,j}:=(-1)^iJ(\omega_{j,i}) \;,
$$
where the subscript $i,j$ denotes the $(i,j)$-th component of a form. 

\bigskip
Extra structures on $\Omega^{\bullet,\bullet}$ are a closed integral
and the corresponding Hodge $*$-operator.

Notice that $\Omega^{2,2}(\CP^2_q)=\Aq$ is a free module of rank one (cf.~Table \ref{tab}),
with basis the element $1$. Let $\mathtt{vol}$ be the form with all components
equal to zero but for the one in degree $4$, which is $1$. We think of this as
the volume form and define an integral by 
\begin{equation}\label{eq:nint}
\nint\omega:=
\inner{\mathtt{vol},\omega}
=\varphi(\omega_{2,2}) \;,\qquad\forall\;\omega\in\Omega^{\bullet,\bullet}(\CP^2_q)\; ; 
\end{equation}
the volume is normalized: $\int\mkern-16mu-\;\mathtt{vol}=1$. Since the differentials $\de$
and $\deb$ are given by the (right) action of elements of $\Uq$ in the kernel of
the counit $\epsilon$, the integral is automatically closed, i.e.
$$
\nint\deb\omega=\nint\de\omega=0 \; ,
$$
a simple consequence of the invariance of the Haar state: $\varphi(a\za x)=\epsilon(x)\varphi(a)$.
The Hodge star operator is the linear operator
$\ast_H:\Omega^{i,j}(\CP^2_q)\to\Omega^{2-j,2-i}(\CP^2_q)$ defined by
\begin{equation}\label{eq:4.3}
\omega^*\!\wprod\omega'=(\ast_H\,\omega,\omega')\texttt{vol} \;,
\end{equation}
for all $\omega\in\Omega^{i,j}(\CP^2_q)$ and $\omega'\in\Omega^{2-j,2-i}(\CP^2_q)$.
The product of forms can be defined in such a way that $\ast_H^2 \omega =(-1)^{\mathrm{dg}(\omega)}\omega$.
{}From this last property, it follows the equality
$$
\dd^\dag= \ast_H{}\dd{}\ast_H \;.
$$
Also, denoting $\mathfrak{e}(\omega)$ the left `exterior product': 
$\mathfrak{e}(\omega)\omega':=\omega\!\wprod\omega'$, and $\mathfrak{i}(\omega)=\mathfrak{e}(\omega)^\dag$
the `contraction' by $\omega$, 
by integrating both sides of \eqref{eq:4.3}, the non degeneracy of the scalar product leads to a nice geometrical interpretation of the Hodge star:
$$
\ast_H\omega=\mathfrak{i}(\omega^*)\mathtt{vol} \;.
$$

\subsection{The spin$^{\lowercase{c}}$ structure}
In analogy to the $q=1$ case, we call $\HH_+$ (resp.~$\HH_-$) the Hilbert
space completion of $\Omega^{0,0}(\CP^2_q)\oplus\Omega^{0,2}(\CP^2_q)$ (resp.~$\Omega^{0,1}(\CP^2_q)$),
and set $\HH:=\HH_+\oplus\HH_-$ with the obvious grading $\gamma:=1\oplus -1$.
The operator $\deb$, as its formal adjoint $\deb^\dag$, are well defined on
the dense subspace $\Omega^{0,\bullet}(\CP^2_q)\subset\HH$ and anticommute with
the grading. Also, if we represent $\Aq$ on $\HH$ via the left multiplication,
the graded Leibniz rule implies that the commutator $[\deb,f]$ can be closed
to a bounded operator on $\HH$ for any $f\in\Aq$. Indeed,
$$
[\deb,f]=\mathfrak{e}(\deb f)
$$
is the operator of `exterior product' by the differential of $f$. Also,
$[\deb^\dag,f]=-[\deb,f^*]^\dag$ tells that the operator $D$ given by
$$
D\omega :=(\deb^\dag v,\deb a+s\deb^\dag b,s\deb v) \;,\qquad
\forall\;\omega=(a,v,b)\in\Omega^{0,\bullet}(\CP^2_q) \;,
$$
has bounded commutators with the algebra $\Aq$, whatever is the value of
$s\in\R$. Notice that $D$ anticommutes with the grading, and it is symmetric.
Self-adjoint extensions of $D$ are in bijections with selfadjoint
extensions of its absolute value, and densely defined positive operators have a canonical
self-adjoint extension, namely the Friedrichs extension.
In order to claim that $(\Aq,\HH,D,\gamma)$ is an even unital
spectral triple (equivariant, since the differential calculus is equivariant)
it remains to check that $D$ has a compact resolvent. This is done in \cite{DDL08b},
for the particular choice $s=\sqrt{[2]_q/2}$, by diagonalizing $D^2$. Indeed,
for $s=\sqrt{[2]_q/2}$ it results
$$
D^2\omega=[2]_q^{-1}\omega(\mathcal{C}_q-2)
$$
for all $\omega\in\Omega^{0,\bullet}(\CP^2_q)$, where $\mathcal{C}_q$ is a suitable element in
the center of $\Uq$, the $q$-analogue of the quadratic Casimir of $\SU(3)$ (see~\cite{DDL08b}).
For central elements left and right canonical actions on the dual Hopf algebra coincide,
so in this case $\omega\za\mathcal{C}_q=\mathcal{C}_q\az\omega$ for all $\omega\in\Omega^{0,\bullet}(\CP^2_q)$.
Now, the decomposition of the \emph{left} $\Uq$-modules $\Omega^{0,k}(\CP^2_q)$ into irreducible
representations is known, as well as the eigenvalue of $\mathcal{C}_q$ in any irreducible
representation. In this way, we get eigenvalues (and their multiplicities) of $D^2$.
If $\lambda\geq 0$ is an eigenvalue of $D^2$, and $\mu$ its multiplicity,
the operator $D$ has eigenvalues $+\lambda$ and $-\lambda$ with multiplicities
$\mu_+$ and $\mu_-$ satisfying $\mu_++\mu_-=\mu$. But the equation $\gamma D=-D\gamma$
means that $\gamma$ sends an eigenvector of $D$ to another eigenvector with opposite
eigenvalue, i.e.~$\mu_+=\mu_-=\mu/2$. With this, the computation of the spectrum of
$D$ is concluded. We give here the result, and refer to \cite{DDL08b} for the details.

\begin{prop}
The operator $D$ has one dimensional kernel, and non-zero
eigenvalues
\begin{align*}
& \pm\sqrt{\tfrac{2}{\,[2]_q\!}[n]_q[n+2]_q}
&& \textrm{\textup{with multiplicity}} \quad (n+1)^3\;,\\
& \pm\sqrt{[n+1]_q[n+2]_q} \hspace{-1cm}
&& \textrm{\textup{with multiplicity}} \quad \tfrac{1}{2}n(n+3)(2n+3) \;,
\end{align*}
for $n\geq 1$.
\end{prop}

\noindent
As a corollary, $(1+D^2)^{-1}$ is compact. In fact,
$(1+D^2)^{-\varepsilon}$ is of trace class for any $\varepsilon>0$.
In this case we usually say that the summability (or the `metric dimension')
of $D$ is $0^+$.

\subsection{Monopoles}
Let $(\Omega^\bullet,\dd)$ be a differential calculus on the algebra $\A$ and $\E$
a right $\A$-module with an $\A$-valued Hermitian structure $(\,,\,)$. A \emph{connection}
on $\E$ (compatible with the Hermitian structure) is a linear map
$\nabla:\E \otimes_{\A}\Omega^\bullet \to\E \otimes_{\A}\Omega^{\bullet+1}$
which satisfies the Leibniz rule
$$
\nabla(\eta \omega)=(\nabla\eta)\omega+(-1)^{\mathrm{dg}(\eta)}\eta(\dd\omega)
$$
and the condition
$$
(\nabla\eta,\xi)+(-1)^{\mathrm{dg}(\eta)}(\eta,\nabla\xi)=\dd(\eta,\xi)
$$
for all $\eta,\xi\in\E\otimes_{\A}\Omega^\bullet$ and $\omega\in\Omega^\bullet$. The Hermitian structure is extended in a natural way:  if
$\eta,\eta'\in\E$, $\omega\in\Omega^{i}$ and $\omega'\in\Omega^{j}$ we define
$$
(\eta\otimes\omega,\eta'\otimes\omega'):=\omega^*(\eta,\eta')\omega'\,\in\Omega^{i+j} \;.
$$
On a finitely generated projective module $\E=e\A^k$ there always exists a connection
$\nabla_e$, called \emph{Grassmannian connection}, and given by
$$
\nabla_e\eta=e\, \dd \eta \;,
$$
where $\dd$ act on $\A^k$ diagonally, and matrix multiplication is understood.
Any other connection differs from $\nabla_e$ by an element in
$\mathrm{Hom}(\E, \E\otimes_{\A}\Omega^1)$.

On $\CP^2_q$, we can use the isomorphism given in previous
sections to transport on $\Sigma_{0,N}$ the Grassmannian connection of $\E=P_N\Aq^{\rPN}$,
with $P_N$ in \eqref{mon-pro}.
For $\eta\in\Sigma_{0,N}\otimes_{\Aq}\Omega^{i,j}(\CP^2_q)$ we get that
$$
\nabla_{\!N}\eta=\Psi_N^\dag\dd(\Psi_N\eta) \;.
$$
The connection $\nabla_{\!N}$ is left $\Uq$-invariant. Its curvature is the operator of left
multiplication by the invariant $2$-form $\nabla_{\!N}^2\in\Omega^{1,1}(\CP^2_q)$ given by
$$
\nabla_{\!N}^2=\Psi_N^\dag(\dd P_N\wprod\dd P_N)\Psi_N \;.
$$
With an explicit computation we find that the curvature is anti-self-dual,
$$
\ast_H\,\nabla_{\!N}^2= -\nabla_{\!N}^2 \;,
$$
and proportional to the curvature of the `tautological' bundle $\Sigma_{0,1}$:
\begin{equation}\label{curvN1}
\nabla_{\!N}^2=q^{N-1} [N] \,\nabla_1^2 \;,
\end{equation}
for any $N\in\Z$. Anti-selfdual connections would be selfdual for the  reversed orientation.

To the connection $\nabla_{\!N}$ there is associated a Laplacian $\Box_N=\nabla_N^\dag\nabla_{\!N}$, acting on 
$\Sigma_{0,N}$, with $\nabla_N^\dag$ the Hermitian conjugate of $\nabla_{\!N}$. As in the case of the
Dirac operator, also the Laplacian can be related to the Casimir $\mathcal{C}_q$. We have that
$$
\Box_N = 
q^{-\frac{3}{2}}\,\frac{q^{\frac{3}{2}}+q^{-\frac{3}{2}}}{q^{\frac{N}{3}}+q^{-\frac{N}{3}}} 
\Big( \mathcal{C}_q-[\tfrac{1}{3}N]_q^2-[\tfrac{1}{3}N+1]_q^2-[\tfrac{2}{3}N+1]_q^2 \Big)
+[2]_q[N]_q \;,
$$
and from this, the spectrum is readily computed.  The spectrum $\{\lambda_{n,N}\}_{n\in\N}$
of $\Box_N$ is given by (we omit the multiplicities)
\begin{align*}
\lambda_{n,N}&=
(1+q^{-3})[n]_q[n+N+2]_q+[2]_q[N]_q &&\mathrm{if}\;N\geq 0\;,\\
\lambda_{n,N}&=
(1+q^{-3})[n+2]_q[n-N]_q+[2]_q[N]_q &&\mathrm{if}\;N\leq 0\;.
\end{align*}
with $n\in\N$. It is worth stressing that the spectrum is not invariant under
the exchange $N \leftrightarrow -N$, not even when sending $q \leftrightarrow q^{-1}$.

\section{Classical and quantum characteristic classes}

A natural question to ask is whether two given finitely generated projective modules are
equivalent. In the presence of symmetries, one can consider equivalence both
in $K$-theory or in equivariant $K$-theory. In the present case the former leads
to integer-valued invariants, discussed in Sec.~\ref{sec:ci},
and the latter to invariants whose values are $q$-numbers, discussed
in Sec.~\ref{sec:qi}.

\subsection{Classical characteristic classes}\label{sec:ci}
Equivalence classes of finitely generated projective (left or right) modules over an algebra $\A$ --
the algebraic counterpart of vector bundles -- are elements of the group $K_0(\A)$.  Equivalence classes
of even Fredholm modules -- the algebraic counterpart of `fundamental classes' -- give a dual group
$K^0(\A)$. There is a map from $K$-theory to cyclic homology, and a map from $K$-homology to cyclic cohomology.
The pullback of the pairing between cyclic cycles and cyclic cocycles gives a pairing
between $K^0(\A)$ and $K_0(\A)$: thus, we can think of an Fredholm
module as a map $K_0(\A)\to\C$. Such a map is actually integer valued,
coinciding with the index of a suitable Fredholm operator~\cite{Con94}.

For the quantum projective plane, $K_0(\Aq)\simeq\Z^3 \simeq K^0(\Aq)$. The result for 
$K$-theory can be proved viewing the corresponding $C^*$-algebra as the Cuntz--Krieger
algebra of a graph \cite{HS02}. The group $K_0$ is given as the cokernel of the incidence
matrix canonically associated with the graph; the result for $K$-homology can be proven
using similar techniques: the groups $K^0$ is now given as the kernel of the transposed
matrix \cite{Cun84}.
Thus finitely generated projective (left or right) modules over $\Aq$ are classified by three
integers. Here we give the three Fredholm modules for $\Aq$ and the corresponding 
invariants $K_0(\Aq)\to\Z$ and refer to \cite{DL09} for
more details.

There is a (unique) non-trivial character $\chi_0:\Aq\to\C$ , defined on
generators by $\chi(p_{ij})=\delta_{i3}\delta_{j3}$. Thinking of $\C$ as functions over
a single point, we interprete $\chi_0$ as (the dual of) the inclusion of a point inside
$\CP^2_q$. We have a first Fredholm module $(\pi_0,\HH_0,F_0)$:
\begin{equation}\label{eq:rank}
\mathrm{ch}^0_{(\pi_0,\HH_0,F_0)}:K_0(\Aq)\to\Z \;,\qquad
\mathrm{ch}^0_{(\pi_0,\HH_0,F_0)}([e])=\tr\,\chi_0(e) \;,
\end{equation}
Its geometrical meaning is the following: the rank of a vector bundle is the dimension of the fiber at any point $x$ of the space, and this coincides with the trace of the corresponding projection evaluated at $x$. With only one `classical point',  the map in \eqref{eq:rank} computes the rank of the restriction of the vector bundle to this classical point. If the module is free, then \eqref{eq:rank} gives its rank.

An irreducible representation $\chi_1$ of $\Aq$ on $\ell^2(\N)$ is given by
\begin{align*}
\chi_1(p_{1i}) &=0 \;,\qquad\forall\;i=1,2,3\,,\\
\chi_1(p_{22})\ket{n} &=q^{2n}\ket{n} \;, \\
\chi_1(p_{23})\ket{n} &=q^{n+1}\sqrt{1-q^{2(n+1)}}\ket{n+1} \;,
\end{align*}
and the operator $\chi_1(p_{ij})-\chi_0(p_{ij})$ is of trace class for all $i,j$.
The associated 
character coming from a Fredholm module $(\pi_1,\HH_1,F_1)$, is
\begin{equation}\label{eq:charge}
\mathrm{ch}^0_{(\pi_1,\HH_1,F_1)}:K_0(\Aq)\to\Z \;,\qquad
\mathrm{ch}^0_{(\pi_1,\HH_1,F_1)}([e])=\tr_{\ell^2(\N)\otimes\C^m}\,(\chi_1-\chi_0)(e) \;,
\end{equation}
where $m$ is the size of the matrix $e$. The value in \eqref{eq:charge} depends only on the restriction of the `vector bundle' to the subspace $\CP^1_q$, and could then be called the \emph{monopole charge} (the 1st Chern number of the bundle).

For the third Fredholm module we take as Hilbert space $\HH_2$ two copies of  the linear span of orthonormal vectors $\ket{\ell,m}$, with $\ell\in\frac{1}{2}\N$ and $\ell+m\in\N$.
The grading $\gamma_2$ and the operator $F_2$ are the obvious ones. It remains to describe the representation
$\pi_2=\pi_+\oplus\pi_-$. Modulo traceclass operators one has:
\begin{align*}
\pi_+(p_{11}) &\sim\pi_+(p_{12})\sim\pi_+(p_{13})\sim 0 \;,\\
\pi_+(p_{22})\ket{\ell,m} &\sim\begin{cases}
q^{2(\ell+m)}\ket{\ell,m}    &\mathrm{if}\;m\leq\ell \;,\\
0                            &\mathrm{if}\;m>\ell \;,
\end{cases} \\
\pi_+(p_{23})\ket{\ell,m} &\sim\begin{cases}
q^{\ell+m+1}\sqrt{1-q^{2(\ell+m+1)}}\ket{\ell,m+1}     &\mathrm{if}\;m\leq\ell-1 \;,\\
0                                                      &\mathrm{if}\;m\geq\ell \;.
\end{cases}
\end{align*}
We define the subrepresentation $\pi_-$ by adding multiplicities to $\chi_1$. On the generators:
\begin{align*}
\pi_-(p_{11}) &=\pi_-(p_{12})=\pi_-(p_{13})=0 \;,\\
\pi_-(p_{22})\ket{\ell,m} &=q^{2(\ell+m)}\ket{\ell,m} \;,\quad
\pi_-(p_{23})\ket{\ell,m} &=q^{\ell+m+1}\sqrt{1-q^{2(\ell+m+1)}}\ket{\ell,m+1} \;.
\end{align*}
On each invariant subspace with a fixed $\ell$, putting $n=\ell+m$ one just recovers  the representation $\chi_1$. Since $\sum_{m>\ell}q^{2(\ell+m)}=(1-q^4)^{-2}$ is finite, on the subspace $m>\ell$ the operators $\pi_-(p_{22})$ and $\pi_-(p_{23})$ are trace class, and so $\pi_+(a)-\pi_-(a)$ is of trace class as well for all $a\in\Aq$: the Fredholm module is $1$-summable. The corresponding character is 
\begin{equation}\label{eq:instnumb}
\mathrm{ch}^0_{(\pi_2,\HH_2,F_2)}:K_0(\Aq)\to\Z \;,\qquad
\mathrm{ch}^0_{(\pi_2,\HH_2,F_2)}([e])=\tr_{\HH_2\otimes\C^m}\,(\pi_+-\pi_-)(e) \;,
\end{equation}
where $m$ is the size of the matrix $e$. The above replaces the 2nd Chern class of the module.

A peculiarity of the quantum case is that in the construction of the characters one needs only to consider the irreducible representations: at $q=1$ irreducible representation are all $1$-dimensional and give only one of the generators of the $K$-homology (the trivial Fredholm module). An additional true `quantum effect' is that the three characters are all traces on $\Aq$: all relevant information leaves in degree zero in cyclic homology. In contrast, for the classical $\CP^2$ one needs to consider homology classes (de Rham currents) in degree $0$, $2$ and $4$.

Pairing the projection $P_N$ in \eqref{mon-pro} with the three Fredholm modules gives, for $N\in\Z$, 
$$
\mathrm{ch}^0_{(\pi_0,\HH_0,F_0)}([P_N])=1 \;,\quad
\mathrm{ch}^0_{(\pi_1,\HH_1,F_1)}([P_N])=N \;,\quad
\mathrm{ch}^0_{(\pi_2,\HH_2,F_2)}([P_N])=\tfrac{1}{2}N(N+1) \;.
$$
We have already mentioned that $K_0(\Aq)\simeq K^0(\Aq)\simeq\Z^3$.
The matrix of the pairings between the three Fredholm modules and
the classes $[1]$, $[P_{-1}]$ and $[P_1]$ is
$$
\begin{pmatrix}
1 & 1 & 1 \\ 0 & -1 & 1 \\ 0 & 0 & 1
\end{pmatrix} ,
$$
which is invertible in $GL(3,\Z)$ with inverse
$$
\begin{pmatrix}
1 & 1 & -2 \\ 0 & -1 & 1 \\ 0 & 0 & 1
\end{pmatrix} .
$$
This proves that the above-mentioned Fredholm modules (resp.~projections) are a basis
of $K^0(\Aq)\simeq\Z^3$ (resp.~$K_0(\Aq)\simeq\Z^3$) as a $\Z$-module, which is equivalent
to saying that they generate them as abelian groups.

\subsection{Quantum characteristic classes}\label{sec:qi}
Classically, invariants of vector bundles are computed by integrating powers of the curvature of
a connection on the bundle, the result being independent of the particular chosen connection.
On the other hand, in order to integrate the curvature of a connection on the quantum projective space 
$\CP^2_q$ one needs `twisted integrals'; the result, as we shall see, is no longer an integer but
rather its $q$-analogue.

We start with some general facts about equivariant algebraic $K$-theory
and $K$-homology and corresponding Chern-Connes characters.
For an homogeneous space, the equivariant topological $K^0$-group is defined as the Grothendieck group
of the abelian monoid whose elements are equivalence classes of equivariant vector bundles. It has an
algebraic version that can be generalized to noncommutative algebras. Let $\U$ be a bialgebra and
$\A$ an $\U$-module algebra. Equivariant vector bundles are replaced by one sided (say left)
$\A\rtimes\U$-modules that are finitely generated and projective as (left) $\A$-modules; these
will be simply called  ``equivariant projective modules''. Any such a module is given by a pair
$(e,\sigma)$, where $e$ is an $k\times k$ idempotent with entries in $\A$, and
$\sigma:\U\to\mathrm{Mat}_k(\C)$ is a representation and the following compatibility requirement
is satisfied (see e.g.~\cite[Sec.~2]{DDL08}):
\begin{equation}\label{eq:cov}
( h_{(1)}\az e)  \sigma(h_{(2)})^t=\sigma(h)^t e \;, \qquad \textup{for all} \quad h\in\U\;,
\end{equation}
with `$\phantom{|}^t$' denoting transposition.
The corresponding module $\E=\A^ke$ is made of elements  $v=(v_1,\ldots,v_k)\in\A^k$ 
in the range of the idempotent, $ve=v$, with module structures
$$
(a.v)_i:=av_i\;,\qquad
(h.v)_i:=\sum\nolimits_{j=1}^k(h_{(1)}\az v_j)\sigma_{ij}(h_{(2)})\;, 
\qquad \textup{for}   \quad a\in\A  \quad \textup{and} \quad  h\in\U \,.
$$
An equivalence between any two equivariant modules is simply an invertible left $\A\rtimes\U$-module map between them.
The group $K_0^{\U}(\A)$ is defined as the Grothendieck group of the abelian monoid whose elements are equivalence
classes of $\U$-equivariant projective $\A$-modules; the monoid operation is the direct sum, as usual.

In the `non-equivariant' case the equivalence relation on finitely generated projective modules
(or vector bundles) is reformulated as an equivalence relation on the corresponding idempotents,
the so-called Murray-von Neumann equivalence. Similarly in the equivariant case, starting with
the equivalence between equivariant projective modules, one is led to the following equivalence
relations on `equivariant' idempotents~\cite{DL09}.

\begin{prop}
Two equivariant projective modules $\A^ke$ and $\A^{k'}e'$ are equivalent if{}f $e=uv$ and $e'=vu$ for some
$u\in\mathrm{Mat}_{k\times k'}(\A)$ and $v\in\mathrm{Mat}_{k'\times k}(\A)$ satisfying the equivariance conditions
$$
(h_{(1)}\az u)\sigma'(h_{(2)})^t=\sigma(h)^tu \;,\qquad
(h_{(1)}\az v)\sigma(h_{(2)})^t=\sigma'(h)^tv \;.
$$
\end{prop}

There is a natural map from equivariant $K$-theory to equivariant cyclic homology given for instance in \cite{NT03}. We adapt the construction there to our situation. One starts with the space $\mathrm{Hom}_{\C}(\U,\A^{n+1})$ 
of $\C$-linear maps from $\U$ to $\A^{n+1}$, 
and for $i=0,\ldots,n$, defines operations $b_{n,i}:\mathrm{Hom}_{\C}(\U,\A^{n+1})\to \mathrm{Hom}_{\C}(\U,\A^n)$ by  
\begin{align*}
b_{n,i}(a_0\otimes a_1\otimes\ldots\otimes a_n)(x)&:=
(a_0\otimes\ldots\otimes a_ia_{i+1}\otimes\ldots\otimes a_n)(x) \;,\qquad\mathrm{if}\;i\neq n\;, \\
b_{n,n}(a_0\otimes a_1\otimes\ldots\otimes a_n)(x)&:=\bigl((x_{(1)}\az a_n)a_0\otimes a_1\otimes\ldots\otimes a_{n-1}\bigr)(x_{(2)}) \;,
\end{align*}
and an operation $\lambda_n:\mathrm{Hom}_{\C}(\U,\A^{n+1})\to \mathrm{Hom}_{\C}(\U,\A^{n+1})$:
$$
\lambda_n(a_0\otimes a_1\otimes\ldots\otimes a_n)(x) :=(-1)^n\bigl((x_{(1)}\az a_n)\otimes a_0\otimes a_1\otimes\ldots\otimes a_{n-1}\bigr)(x_{(2)}) \;.
$$
The maps $b_{n,i}$ make up a presimplicial module --  one checks that $b_{n-1,i}b_{n,j}=b_{n-1,j-1}b_{n,i}$ for all $0\leq i<j\leq n$ --, so that 
$$
b_n:=\sum\nolimits_{i=0}^n(-1)^ib_{n,i}
$$
is a boundary operator \cite{Lod97}.
The Hopf algebra $\U$ acts on $\A^{n+1}$ via the rule
$$
h\aaz (a_0\otimes a_1\otimes\ldots\otimes a_n):=h_{(1)}\az a_0\otimes h_{(2)}\az a_1\otimes\ldots\otimes h_{(n+1)}\az a_n \;,
$$
and $C^{\U}_n(\A)$ will denote the collection of elements $\omega\in\mathrm{Hom}_{\C}(\U,\A^{n+1})$ which are `equivariant', meaning that
$$
(h_{(1)}\aaz\omega)(xh_{(2)})=\omega(hx) \;,
$$
for all $h,x\in\U$. The operators $b_{n,i}$ commute with the action of $\U$, and it makes sense to consider the complex of equivariant maps. 
The cyclic operator $\lambda_n$ commutes with the action of $\U$, thus it descends to an operator on $C^{\U}_n(\A)$ as well.
Finally, with 
$$
b'_n:=\sum\nolimits_{i=0}^{n-1}(-1)^nb_{n,i} \,,
$$  
it holds that $b_n(1-\lambda_n)=(1-\lambda_{n-1})b'_n$, which says that the boundary operator $b_n$ maps the space $C^{\U}_n(\A)/\mathrm{Im}(1-\lambda_n)$ into $C^{\U}_{n-1}(\A)/\mathrm{Im}(1-\lambda_{n-1})$.
The homology of this last complex is called ``$\U$-equivariant cyclic homology'' of $\A$, with corresponding homology groups $H\!C^{\U}_n(\A)$.

Next, for $\sigma:\U\to\mathrm{Mat}_k(\C)$ a representation as in \eqref{eq:cov} above, consider the set
$$
\mathrm{Mat}_k^\sigma(\A):=\big\{a\in\mathrm{Mat}_k(\A)\; \big|\; 
\left(h_{(1)}\az a \right) \sigma(h_{(2)})^t=\sigma(h)^t a\;, \;\; \forall\;h\in\U \,\big\} \;.
$$
This is a subalgebra of $\mathrm{Mat}_k(\A)$; given indeed any two of its elements $a,b$ one has: 
$$
\left(h_{(1)}\az (ab)\right) \sigma(h_{(2)})^t=
\left((h_{(1)}\az a)(h_{(2)}\az b) \right) \sigma(h_{(3)})^t \\ =
(h_{(1)}\az a) \, \sigma(h_{(2)})^t \, b=\sigma(h)^t \, ab\;.
$$
Moreover, $\sigma$-equivariant $k\times k$ idempotents as in \eqref{eq:cov} are elements of $\mathrm{Mat}^\sigma_k(\A)$.
Due to the definition of $\mathrm{Mat}^\sigma_k(\A)$ there exists a map $\tr_\sigma:\mathrm{Mat}^\sigma_k(\A)^{n+1}\to C^{\U}_n(\A)$ given by
\begin{align*}
\tr_\sigma(a_0\otimes a_1\otimes\ldots\otimes a_n)(x)
&:=\tr\bigl(a_0\dotimes a_1\dotimes\ldots\dotimes a_n\sigma(x)^t\bigr) \\
&=\sum_{i_0,i_1,\ldots,i_{n+1}}(a_0)_{i_0i_1}\otimes (a_1)_{i_1i_2}\otimes\ldots\otimes
(a_n)_{i_ni_{n+1}}\sigma(x)_{i_0i_{n+1}} \;,
\end{align*}
where $\dotimes$ denotes composition of the tensor product over $\C$ with matrix multiplication.
The map $\tr_\sigma$ is a morphism of differential complexes, mapping the complex of the cyclic homology of
$\mathrm{Mat}^\sigma_k(\A)$ to the complex of the $\U$-equivariant cyclic homology of $\A$. This construction
is completely analogous to the ``non-equivariant'' case, cf.~\cite[Cor.~1.2.3]{Lod97}.
At this point, one can repeat verbatim the proof of Thm~8.3.2 in \cite{Lod97}, replacing the ring
$R:=\mathrm{Mat}_k(\A)$ there, with $\mathrm{Mat}^\sigma_k(\A)$ (which is still a ring) and replacing
the generalized trace map there, with $\tr_\sigma$, to prove the following theorem.

\begin{thm}
A map $\mathrm{ch}^n:K_0^{\U}(\A)\to H\!C^{\U}_n(\A)$ is defined by
$$
\mathrm{ch}^n(e,\sigma):=\tr_\sigma(e^{\otimes n+1}) \;.
$$
\end{thm}

\noindent
We remark that what we denote here $H\!C_n$ and call cyclic homology is Connes' first version
of cyclic homology, i.e.~the homology of Connes' complex denoted $H^\lambda_n$ in \cite{Lod97}.

Modulo a normalization, the map $\mathrm{ch}^n(e,\sigma)(1)$ is the usual Chern-Connes
character in cyclic homology (and no $\sigma$'s in the formul{\ae}). 
In general, one fixes a group-like element $K\in\U$ calling $\eta$
the corresponding automorphism of $\A$, $\eta(a):=K\az a$ for all $a\in\A$.
Next, we call ${_\eta}\A$ the $\A$ bimodule that is $\A$ itself as a vector
space, but has a left module structure `twisted' with $\eta$:
$$
a.m=\eta(a)m\qquad\mathrm{and}\qquad m.a=ma
$$
for all $a\in\A$ and $m\in{_\eta}\A$, where the dot denotes the bimodule structure
while without dot we mean the product in $\A$. The groups $H\!C^{\U}_\bullet(\A)$ can
be paired with the Hochschild cohomology group of $\A$ with coefficients in ${_\eta}\A$.
Indeed, the pairing 
\begin{equation}\label{twdprng}
\inner{\,,\,}:\mathrm{Hom}_{\C}(\A^{n+1},\C)\times\mathrm{Hom}_{\C}(\U,\A^{n+1})\to\C \, , \qquad 
\inner{\tau,\omega}:=\tau\bigl(\omega(K)\bigr) \;,
\end{equation}
when used to compute the dual $b^*_{n,i}:\mathrm{Hom}_{\C}(\A^n,\C)\to\mathrm{Hom}_{\C}(\A^{n+1},\C)$
of the face operators $b_{n,i}$, yields the formul{\ae} 
\begin{align*}
b_{n,i}^*\tau(a_0,a_1,\ldots,a_n)&=\tau(a_0,\ldots,a_ia_{i+1},\ldots,a_n) \;,\qquad\mathrm{if}\;i\neq n\;,\\
b_{n,n}^*\tau(a_0,a_1,\ldots,a_n)&:=\tau\bigl(\eta(a_n)a_0,a_1,\ldots,a_{n-1}\bigr) \;.
\end{align*}
The above are just the face operators of the Hochschild cohomology $H^\bullet(\A,{_\eta}\A)$ of $\A$ with
coefficients in ${_\eta}\A$ (cf.~\cite{Lod97}). Thus, the pairing in \eqref{twdprng} descends to a pairing 
$$
H^n(\A,{_\eta}\A)\times H\!C^{\U}_n(\A)\to\C \,. 
$$
The maps $\mathrm{ch}^n(e,\sigma)((K_1 K_2)^{-4})$ are what we shall use for $\CP^2_q$ in the following section.

\subsection{The case of $\CP^2_{\lowercase{q}}$}
As mentioned, for $\CP^2_q$ we set $K=(K_1K_2)^{-4}$ the element implementing the modular automorphism.
Indeed, the Haar state of $\Oq$ satisfies (cf.~\cite{KS97})
$$
\varphi(ab)=\varphi\bigl( (K\az b\za K ) a\bigr) \;,\qquad \textup{for} \quad a,b\in\Oq\;,
$$
that when $a,b\in\Aq$ means
$$
\varphi(ab)=\varphi\bigl( (K\az b) a\bigr)=\varphi\bigl( \eta(b) a\bigr) \;.
$$
That is, the restriction of the Haar state to $\Aq$ is the representative of a class in the cohomology 
$H^0(\Aq,{_\eta}\Aq)$.
On the other hand, with the integral defined in \eqref{eq:nint} by using the Haar state as well, an element
$[\tau_4]\in H^4(\Aq,{_\eta}\Aq)$ is constructed as 
$$
\tau_4(a_0,\ldots,a_4):=\nint a_0\dd a_1\wprod\ldots\wprod\dd a_4 \;.
$$

\medskip
A $2$-cocycle can be defined in a similar way. Recall that elements of $\Omega^{1,1}(\CP^2_q)$
have the form $\omega=(\alpha,\alpha_4)$, with $\alpha_4\in\Aq$. Let $\pi:\Omega^{1,1}(\CP^2_q)\to\Aq$ be the projection onto the second component $\pi(\omega)=\alpha_4$, and extend it to a projection $\pi:\Omega^2(\CP^2_q)\to\Aq$
by setting $\pi(\omega)=0$ if $\omega\in\Omega^{0,2}$ or $\omega\in\Omega^{2,0}$. 
The map $\pi$ is an $\Aq$-bimodule map. Then, the map
$$
\tau_2(a_0,a_1,a_2):=\varphi\circ\pi(a_0\dd a_1\wprod \dd a_2)
$$
is the representative of a class $[\tau_2]\in H^2(\Aq,{_\eta}\Aq)$. Indeed, by 
the Leibniz rule,
$$
b^*_3\tau_2(a_0,a_1,a_2,a_3)=\varphi\circ\pi\bigl(
a_0(\dd a_1\wprod \dd a_2)a_3-\eta(a_3)a_0(\dd a_1\wprod \dd a_2)\bigr) \;.
$$
Being $\pi$ a bimodule map we get in turn 
$$
b^*_3\tau_2(a_0,a_1,a_2,a_3)=\varphi\bigl(
a_0\pi(\dd a_1\wprod \dd a_2)a_3-\eta(a_3)a_0\pi(\dd a_1\wprod \dd a_2)\bigr) \;,
$$
which is zero by the modular property of the Haar state. 

Both classes $[\tau_4]$ and $[\tau_2]$ will be proven to be not trivial by pairing them with the 
monopole projections \eqref{mon-pro}. Firstly, we observe that 
the projection $P_N=\Psi_N\Psi_N^\dag$ in \eqref{mon-pro} is equivariant with respect to the
representation $\sigma^N$ given by 
$$
\sigma^N(h):= \begin{cases}
\rho^{0,N}(S(h))^t &\; \mathrm{if}\;N\geq 0\;,\\
\sigma^N(h):=\rho^{-N,0}(S(h))^t &\; \mathrm{if}\;N<0\;. 
\end{cases}
$$
The pairing of $[\tau_4]$ with $\mathrm{ch}^4(P_N,\sigma^N)$ is
$$
\inner{\tau_4,\mathrm{ch}^4(P_N,\sigma^N)}=\nint\tr\bigl(P_N(\dd P_N)^4\sigma^N(K_1^{-4}K_2^{-4})^t\bigr) \;,
$$
and using the modular properties of the Haar state can be rewritten as
the integral of the square of the curvature
$$
\inner{\tau_4,\mathrm{ch}^4(P_N,\sigma^N)}=q^{-2N}\nint \nabla_N^2\wprod\nabla_N^2 \;.
$$
{}From \eqref{curvN1} we know that $\nabla_N^2=q^{N-1} [N] \,\nabla_1^2$. Thus the value of the corresponding quantum Chern number is proportional to $[N]^2$: 
$$
\inner{\tau_4,\mathrm{ch}^4(P_N,\sigma^N)} = 
\left(q^{-2} \, \nint \nabla_1^2\wprod\nabla^2_1 \right) [N]^2  \;.
$$
For $q=1$, the integral of the square of the curvature is (modulo a global normalization constant)
the \emph{instanton number} of the bundle.

The pairing of $[\tau_2]$ with $\mathrm{ch}^2(P_N,\sigma^N)$ gives
\begin{align*}
\inner{\tau_2,\mathrm{ch}^2(P_N,\sigma^N)}
 &=\varphi\,\tr\bigl(P_N\pi(\dd P_N\wprod\dd P_N)\sigma^N(K_1^{-4}K_2^{-4})^t\bigr) \\
 &=q^{-2N}\varphi\bigl(\Psi_N^\dag\pi(\dd P_N\wprod\dd P_N)\Psi_N\bigr) \;.
\end{align*}
Since $\Psi_N$ are `functions' on the total space of the bundle, we cannot move them
inside $\pi$ (which is an $\Aq$-bimodule map, not an $\Oq$-bimodule map). Nevertheless,
 -- with a little abuse of notations -- the form $\nabla_N^2=\pi(\nabla_N^2)$ is a constant, and
$$
P_N\dd P_N\wprod\dd P_N=\Psi_N\nabla^2_N\Psi_N^\dag
=\pi(\nabla_N^2)\Psi_N\Psi_N^\dag=\pi(\nabla_N^2)P_N \;.
$$
With this, and using $\Psi_N^\dag  P_N\Psi_N=1$, we come to the final formula
$$
\inner{\tau_2,\mathrm{ch}^2(P_N,\sigma^N)}
 =q^{-2N}\varphi\circ\pi(\nabla_N^2) \;.
$$
In \eqref{curvN1}, we have already shown  that $\nabla_N^2=q^{N-1} [N] \,\nabla_1^2$ 
with $ \nabla_1^2=(0,w_1)$. Thus the corresponding quantum Chern number is proportional to 
$q^{-N}[N]$:
$$
\inner{\tau_2,\mathrm{ch}^2(P_N,\sigma^N)}
 = \left( \varphi\circ\pi(\nabla_1^2) \right) q^{-N-1}[N] \,. 
$$ 
At $q=1$ the integral of the curvature is (modulo a global normalization constant)
the \emph{monopole number} of the bundle; it is the same as the first Chern number.

We conclude with a remark. If we pair $\mathrm{ch}^0(P_N,\sigma^N)$ with (the restriction of the Haar state to $\Aq$) 
$\varphi$
we get $\inner{\varphi,\mathrm{ch}^0(P_N,\sigma^N)}=q^{-2N}$.
If $q$ is trascendental, this means that all $[P_N]$ are independent, i.e.~the equivariant
$K_0$-group is infinite dimensional. Indeed, were the classes $[P_N]$ not independent, there
would exist a sequence $\{k_N\}$ of integers -- all zero but for finitely many -- such that
$\sum_Nk_Nq^{-2N}=0$, and $q^{-1}$ would be the root of a non-zero polynomial with integer
coefficients. This is an instance of the general fact (cf. \cite[Thm.~$3.6$]{NT03}) 
that the equivariant $K_0$ group is a free abelian group and generators are, for the present case, 
in bijection with equivalence classes of irreducible corepresentations of $\Kq$.  


\bigskip\medskip
\begin{center}
\textsc{Acknowledgments}.
\end{center}
We are grateful to the organizers of the 
5th ECM Satellite Conference on ``Noncommutative Structures in Mathematics and Physics'',
held at the Royal Flemish Academy, Brussels (Belgium), 22-26 July 2008, for the nice invitations.
GL was partially supported by the `Italian project Cofin06 - Noncommutative geometry, quantum groups and applications'.

\medskip


\providecommand{\bysame}{\leavevmode\hbox to3em{\hrulefill}\thinspace}

\end{document}